\documentclass[10pt]{article}
 \font\smallit=cmti10

\usepackage{amssymb,latexsym,amsmath,epsfig,amsthm}

\makeatletter

\renewcommand\section{\@startsection {section}{1}{\z@}
 {-30pt \@plus -1ex \@minus -.2ex}
 {2.3ex \@plus.2ex}
 {\normalfont\normalsize\bfseries}}

\renewcommand\subsection{\@startsection{subsection}{2}{\z@}
 {-3.25ex\@plus -1ex \@minus -.2ex}
 {1.5ex \@plus .2ex}
 {\normalfont\normalsize\bfseries}}

\renewcommand{\@seccntformat}[1]{\csname the#1\endcsname. }

\makeatother

\newtheorem{theorem}{Theorem}
 \newtheorem{lemma}{Lemma}

\begin{document}

\begin{center}
 \uppercase{\bf Sampling Goldbach Numbers at Random}
 \vskip 20pt
 {\bf Ljuben Mutafchiev}\\
 {\smallit Department of Mathematics and Science, American University in Bulgaria and
 Institute of Mathematics and Informatics of the Bulgarian Academy of Sciences, Blagoevgrad, 2700, Bulgaria}\\
 {\tt ljuben@aubg.edu}\\
 \end{center}
 \vskip 30pt

\centerline{\bf Abstract}

\noindent
 Let $\Sigma_{2n}$ be the set of all partitions of the even
 integers from the interval $(4,2n]$, $n>2,$ into two odd prime
 parts. We select a partition from the
 set $\Sigma_{2n}$ uniformly at random. Let $2G_n$ be the number partitioned by this
 selection. $2G_n$ is sometimes called a Goldbach number. In [6]
 we showed that $G_n/n$ converges weakly to the maximum $T$ of two
 random variables which are independent copies of a unformly
 distributed random variable in the interval $(0,1)$. In this note
 we show that the mean and the variance of $G_n/n$ tend to the
 mean $\mu_T=2/3$ and variance $\sigma_T^2=1/18$ of $T$,
 respectively. Our method of proof is based on generating
 functions and on a Tauberian theorem due to
 Hardy-Littlewood-Karamata.

\pagestyle{myheadings}
 \markright{}
 \thispagestyle{empty}
 \baselineskip=12.875pt
 \vskip 30pt

\section{Introduction and Statement of the Main Result}

A Goldbach number is an even positive integer that can be
represented as the sum of two primes without regard to order. The
representation itself is called a Goldbach partition of the
underlying even positive integer. Let
$\mathcal{P}=\{p_1,p_2,...\}$ be the sequence of all odd primes
arranged in increasing order. For any integer $k>2$, by $Q_2(2k)$
we denote the number of Goldbach partitions of the number $2k$. In
1742 Goldbach conjectured that $Q_2(2k)\ge 1$. This problem
remains still unsolved (for more details, see e.g. [4], Section
2.8 and p. 594, [8], Section 4.6 and [9], Chapter VI). Let
$\Sigma_{2n}$ be the set of all Goldbach partitions of the even
integers from the interval $(4,2n]$, $n>2$. The cardinality of
this set is obviously
\begin{equation}\label{cardinal}
\mid\Sigma_{2n}\mid=\sum_{2<k\le n} Q_2(2k).
\end{equation}
We recently established in [6] that
\begin{equation}\label{asympcard}
\mid\Sigma_{2n}\mid\sim\frac{2n^2}{\log^2{n}}, \quad n\to\infty.
\end{equation}

Consider now a random experiment. Suppose that we select a
partition uniformly at random from the set $\Sigma_{2n}$, i.e. we
assign the probability $1/\mid\Sigma_{2n}\mid$ to each Goldbach
partition of an even integer from the interval $(4,2n]$. We denote
by $\mathbb{P}$ the uniform probability measure on $\Sigma_{2n}$.
Let $2G_n\in (4,2n]$ be the Goldbach number that is partitioned by
this random selection. Using (\ref{asympcard}), it is not
difficult to show that $G_n/n$ converges weakly, as $n\to\infty$,
to a random variable whose cumulative distribution function is
\begin{equation}\label{cdf}
F(u) =\left\{\begin{array}{ll} 0 & \qquad  \mbox {if} \qquad u\le
 0, \\
 u^2 & \qquad \mbox {if}\qquad 0<u<1, \\
 1 & \qquad \mbox {if} \qquad u\ge 1
 \end{array}\right.
\end{equation}
 (see [6], Theorem 2). It can be easily seen that $F(u)$ is the
 distribution function of the random variable
 $T=\max{\{U_1,U_2\}}$, where $U_1$ and $U_2$ are two independent
 copies of a uniformly distributed random variable in the interval
 $(0,1)$.

 The goal of this present note is to determine asymptotically the
 first two moments of a random Goldbach number $2G_n$. We state
 our main result in terms of the expectation $\mathbb{E}$ and
 variance $\mathbb{V}ar$ of the randaom variable $G_n$, both taken
 with respect to the probability measure $\mathbb{P}$. By $\mu_T$
 and $\sigma_T^2$ we denote the expected value and the variance of
 the random variable $T$, respectively. (Its cumulative
 distribution function is given by (\ref{cdf}).)

 \begin{theorem} We have
 $$
 (i) \quad \lim_{n\to\infty}\frac{1}{n}\mathbb{E}(G_n) =\mu_T
 =\frac{2}{3};
 $$

 $$
 (ii) \quad \lim_{n\to\infty}\frac{1}{n^2}\mathbb{V}ar(G_n)
 =\sigma_T^2=\frac{1}{18}.
 $$
 \end{theorem}

 Our method of proof is similar to that in [6]. It is essentially
 based on a generating function identity and on a classical
 Tauberian theorem due to Hardy, Littlewood and Karamata (see e.g.
 [7], Theorem 8.7).

 {\it Remark 1.} Our method of proof and further but similar
 computations yield also asymptotic expressions for moments of
 higher order of the variable $G_n/n$.

 {\it Remark 2.} The proof of Hardy-Littlewood-Karamata Tauberian
 theorem may be found e.g. in [3], Chapter 7.

 Our paper is organized as follows. Section 2 contains some
 preliminaries. The proof of Theorem 1 is given in Section 3.

 \section{Preliminary Results}

 By a prime partition of the positive integer $n$, we mean a way
 of writing it as a sum of primes from the set $\mathcal{P}$
 without regard to order; the summands are called parts. Clearly,
 Goldbach partitions are prime partitions in which the number of
 parts is $2$. Consider now the number $Q_m(n)$ of prime
 partitions of $n$ into $m$ parts ($1\le m\le n$). The bivariate
 generating function of the numbers $Q_m(n)$ is of Euler's type,
 namely,
 \begin{equation}\label{euler}
 1+\sum_{n=1}^\infty z^n\sum_{m=1}^n Q_m(n)x^m
 =\prod_{p_k\in\mathcal{P}} (1-xz^{p_k})^{-1}.
 \end{equation}
 (the proof may be found in [1], Section 2.1). It is clear that,
 for $n>4$, $Q_2(n)$ counts the number of Goldbach partitions of
 $n$ and that $Q_2(n)=0$ if $n$ is odd.

 For any real variable $z$ with $\mid z\mid<1$, we also set
 \begin{equation}\label{fz}
 f(z)=\sum_{p_k\in\mathcal{P}} z^{p_k}.
 \end{equation}
 Differentiating the left- and right-hand sides of (\ref{euler})
 twice with respect to $z$ and setting then $z=0$ and $m=2$, we
 obtain the following identity.

 \begin{lemma} (See [6], Lemma 3.) For $\mid z\mid<1$, we have
 \begin{equation}\label{partition}
 2\sum_{k>2} Q_2(2k)z^{2k}=f^2(z)+f(z^2),
 \end{equation}
 where $f(z)$ is the function defined by (\ref{fz}).
 \end{lemma}
 Further, we will use a Tauberian theorem by
 Hardy-Littlewood-Karamata [3], Chapter 7. We use it in the form
 given by Odlyzko [7], Section 8.2.

 {\bf Hardy-Littlewood-Karamata Theorem.} {\it (See [7; Theorem
8.7, p. 1225].) Suppose that $a_k\ge 0$ for all $k$, and that
$$
g(x)=\sum_{k=0}^\infty a_k x^k
$$
converges for $0\le x<r$. If there is a $\rho>0$ and a function
$L(t)$ that varies slowly at infinity such that
\begin{equation}\label{funcsim}
g(x)\sim (r-x)^{-\rho}L\left(\frac{1}{r-x}\right), \quad x\to r^-,
\end{equation}
then
\begin{equation}\label{sumsim}
\sum_{k=0}^n a_k r^k\sim\left(\frac{n}{r}\right)^\rho
\frac{L(n)}{\Gamma(\rho+1)}, \quad n\to\infty.
\end{equation}}

{\it Remark.} A function $L(t)$ varies slowly at infinity if, for
every $u>0$, $L(ut)\sim L(t)$ as $t\to\infty$.

\section{Proof of the Main Result}

{\it Proof of Theorem 1(i).} First, we notice that the definition
of the random variable $G_n$ implies that
\begin{equation}\label{expect}
\frac{1}{n}\mathbb{E}(G_n)=\sum_{2<k\le n}
\left(\frac{k}{n}\right)\mathbb{P}(G_n=k) =\sum_{2<k\le n}
\left(\frac{2k}{2n}\right)\frac{Q_2(2k)}{\sum_{2<j\le n} Q_2(2j)}.
\end{equation}
The asymptotic behavior of the denominator in the right-hand side
of (\ref{expect}) is completely determined by (\ref{cardinal}) and
(\ref{asympcard}). Differentiating both sides of
(\ref{partition}), we will show next that the series
\begin{equation}\label{firstder}
\sum_{k>2} (2k)Q_2(2k)z^{2k-1} =f(z)f^\prime(z)+zf^\prime(z^2)
\end{equation}
and
\begin{equation}\label{secder}
\sum_{k>2} (2k)(2k-1)Q_2(2k)z^{2k-2} =f^{\prime
2}(z)+f(z)f^{\prime\prime}(z)+f(z^2)+2z^2f^\prime(z^2)
\end{equation}
satisfy the conditions of Hardy-Littlewood-Karamata theorem. We
need the following lemma.

\begin{lemma} Let $f(z)$ be the series defined by (\ref{fz}). Then,
as $z\to 1^-$,
\begin{equation}\label{ef}
f(z)\sim\frac{1}{(1-z)\log{\frac{1}{1-z}}},
\end{equation}
\begin{equation}\label{efprime}
f^\prime(z)\sim\frac{2}{(1-z)^2\log{\frac{1}{1-z}}},
\end{equation}
\begin{equation}\label{efsecond}
f^{\prime\prime}(z)\sim\frac{2}{(1-z)^3\log{\frac{1}{1-z}}}.
\end{equation}
\end{lemma}

{\it Proof.} The proof of (\ref{ef}) is given in [6], Lemma 4.
Here we present a complete proof of (\ref{efprime}). The proof of
(\ref{efsecond}) is similar.

As usual, by $\pi(y)$ we denote the number of primes which do not
exceed the positive number $y$. In (\ref{fz}) we set $z=e^{-t},
t>0$, and apply an argument similar to that given by Stong [10]
(see also [2]). We have
\begin{equation}\label{eminust}
f^\prime(z)\mid_{z=e^{-t}} =\sum_{p_k\in\mathcal{P}} p_k
z^{p_k-1}\mid_{z=e^{-t}} =\int_0^\infty ye^{-yt}d\pi(y)
=I(t)-f(e^{-t}),
\end{equation}
where
\begin{equation}\label{it}
I(t)=t\int_0^\infty ye^{yt}\pi(y)dy.
\end{equation}
In [6] we established that
\begin{equation}\label{efeminust}
f(e^{-t})\sim\frac{1}{t\log{\frac{1}{t}}}, \quad t\to 0^+.
\end{equation}
To find the asymptotic of $I(t)$, we set in (\ref{it}) $y=s/t$ and
obtain
\begin{equation}\label{ittwo}
I(t)=\frac{1}{t}\int_0^\infty se^{-s}\pi(s/t)ds
=\frac{1}{t}(I_1(t)+I_2(t)),
\end{equation}
where
$$
I_1(t)=\int_0^{t^{1/2}} se^{-s}\pi(s/t)ds,
$$

$$
I_2(t)=\int_{t^{1/2}}^\infty se^{-s}\pi(s/t)ds.
$$
For $I_1(t)$ we use the bound $\pi(s/t)\le s/t$. Hence, for enough
small $t>0$, we have
\begin{eqnarray}\label{ione}
& & 0\le I_1(t)\le\frac{1}{t}\int_0^{t^{1/2}} s^2 e^{-s}ds
=\frac{1}{t}\left(-s^2 e^{-s}\mid_0^{t^{1/2}} +2\int_0^{t^{1/2}}
se^{-s}ds\right) \nonumber \\
& & =-e^{-t^{1/2}}+O(t^{-1/2}) =O(t^{-1/2}).
\end{eqnarray}

The estimate of $I_2(t)$ follows from the Prime Number Theorem
with an error term given by Ingham [5], Theorem 23, p.65. Thus,
for $y>1$, we have
$$
\pi(y) =\frac{y}{\log{y}} +O\left(\frac{y}{\log^2{y}}\right).
$$
Hence, for $s\ge t^{1/2}$, we have $\log{s}\ge
-\frac{1}{2}\log{\frac{1}{t}}$ and therefore,
\begin{eqnarray}
& & \pi(s/t)
=\left(\frac{s}{t}\right)\frac{1}{\log{\frac{1}{t}}+\log{s}}
+O\left(\frac{s}{t\left(\log{\frac{1}{t}}+\log{s}\right)^2}\right)
\nonumber \\
& & =\frac{s}{t\log{\frac{1}{t}}}
\left(1+O\left(\frac{\mid\log{s}\mid}{\log{\frac{1}{t}}}\right)\right)
+O\left(\frac{s}{t\log{\frac{1}{t}}}\right) \nonumber \\
& &=\frac{s}{t\log{\frac{1}{t}}}
+O\left(\frac{s(1+\mid\log{s}\mid)}{t\log^2{\frac{1}{t}}}\right)
\nonumber.
\end{eqnarray}
Consequently,
\begin{eqnarray}\label{itwo}
& & I_2(t)=\frac{1}{t\log{\frac{1}{t}}} \int_{t^{1/2}}^\infty s^2
e^{-s}ds +
O\left(\frac{1}{t\log^2{\frac{1}{t}}}\int_{t^{1/2}}^\infty
s^2(1+\mid\log{s}\mid)e^{-s}ds\right) \nonumber \\
& & =\frac{1}{t\log{\frac{1}{t}}}\left(\int_0^\infty s^2 e^{-s}ds
-\int_0^{t^{1/2}} s^2 e^{-s}ds\right)
+O\left(\frac{1}{t\log^2{\frac{1}{t}}}\right) \nonumber \\
& & =\frac{1}{t\log{\frac{1}{t}}}(2-O(t^{1/2}))
+O\left(\frac{1}{t\log^2{\frac{1}{t}}}\right) \nonumber \\
 & & =\frac{2}{t\log{\frac{1}{t}}}
+O\left(\frac{1}{t^{1/2}\log{\frac{1}{t}}}\right)
+O\left(\frac{1}{t\log^2{\frac{1}{t}}}\right)
\sim\frac{2}{t\log{\frac{1}{t}}}, \quad t\to 0^+.
\end{eqnarray}
Combining (\ref{ittwo})-(\ref{itwo}), we get
$$
I(t)\sim\frac{2}{t^2\log{\frac{1}{t}}}, \quad t\to 0^+.
$$
Replacing this asymptotic equivalence and (\ref{efeminust}) into
(\ref{eminust}), we conclude that
$$
f^\prime (z)\mid_{z=e^{-t}}
=I(t)+O\left(\frac{1}{t\log{\frac{1}{t}}}\right)
\sim\frac{2}{t^2\log{\frac{1}{t}}}, \quad t\to 0^+.
$$
Returning to the variable $z$ by $t=\log{\frac{1}{z}}$ ($z<1$), we
obtain
$$
f^\prime(z)\sim -\frac{2}{\left(\log^2{\frac{1}{z}}\right)
\log{\log{\frac{1}{z}}}}, \quad z\to 1^-.
$$
Now, (\ref{efprime}) follows from the obvious observation that
$$
\log{\frac{1}{z}}=-\log{z}=-\log{(1-(1-z))}\sim 1-z, \quad z\to
1^-.
$$

The asymptotic equivalence (\ref{efsecond}) can be obtained in the
same way using the representation
$$
f^{\prime\prime}(z)\mid_{z=e^{-t}} =\int_0^\infty
y(y-1)e^{-yt}d\pi(y).
$$
One can show that this last integral is
$\sim\frac{2}{t^3\log{\frac{1}{t}}}$ as $t\to 0^+$ and then
substitute again $t$ by $-\log{z}$. We omit further details.
$\rule{2mm}{2mm}$

To complete the proof of Theorem 1(i) we recall (\ref{firstder})
and observe that (\ref{efprime}) implies that
$$
f^\prime(z^2) \sim\frac{2}{(1-z^2)\log{\frac{1}{1-z^2}}}
=O\left(\frac{1}{(1-z)^2\log{\frac{1}{1-z}}}\right), \quad z\to
1^-.
$$
So, the main contribution to the asymptotic of the right-hand side
of (\ref{firstder}) is given by the product $f(z)f^\prime(z)$.
Hence, from (\ref{ef}) and (\ref{efprime}) we find that
\begin{eqnarray}\label{asexpect}
& & \sum_{k>2} (2k)Q_2(2k)z^{2k-1}
=\frac{2}{(1-z)^3\log^2{\frac{1}{1-z}}}
+O\left(\frac{1}{(1-z)^2\log{\frac{1}{1-z}}}\right) \nonumber \\
& & \sim\frac{2}{(1-z)^3\log^2{\frac{1}{1-z}}}, \quad z\to 1^-,
\end{eqnarray}
which implies that the series $\sum_{k>2} (2k)Q_2(2k)z^{2k-1}$
satisfies condition (\ref{funcsim}) of Hardy-Littlewood-Karamata
Tauberian theorem with $r=1, \rho=3$ and
$L(t)=\frac{2}{\log^2{t}}$. By (\ref{sumsim}) we obtain
\begin{equation}\label{assum}
\sum_{2<k\le n} (2k)Q_2(2k)
\sim\left(\frac{8}{3}\right)\frac{n^3}{\log^2{n}}, \quad
n\to\infty.
\end{equation}
Furthermore, (\ref{cardinal}) and (\ref{asympcard}) imply that
\begin{equation}\label{initialsum}
(2n)\sum_{2<k\le n} Q_2(2k)\sim\frac{4n^3}{\log^2{n}}, \quad
n\to\infty.
\end{equation}
Dividing (\ref{assum}) by (\ref{initialsum}) we see that the limit
of the right-hand side of (\ref{expect}) is $2/3$ as $n\to\infty$.
$\rule{2mm}{2mm}$

{\it Remark.} Theorem 1(i) also presents a solution of the
following curious problem. Consider a sampling procedure that
consists of two steps. In the first step we select, as previously,
a Goldbach partition from the set $\Sigma_{2n}$, and, in the
second step, we select an even number $2R_n$ from the interval
$(4,2n]$. It is then easy to see that the middle part of
(\ref{expect}) represents the probability $Pr(R_n\le G_n)$ by the
total probability formula. Therefore, Theorem 1(i) implies that
$$
lim_{n\to\infty} Pr(R_n\le G_n)=\frac{2}{3}.
$$

{\it Proof of Theorem 1(ii).} In the same way, we establish that
the leading term in the asymptotic of the right-hand side of
(\ref{secder}) is given by the first two terms $f^{\prime
2}(z)+f(z)f^{\prime\prime}(z)$. By (\ref{ef})-(\ref{efsecond}), we
have
\begin{equation}\label{mainterm}
f^{\prime 2}(z)+f(z)f^{\prime\prime}(z) \sim\frac{6}{(1-z)^4
\log^2{\frac{1}{1-z}}}, \quad z\to 1^-.
\end{equation}
Furthermore, by (\ref{asexpect}),
\begin{eqnarray}
& & \sum_{k>2} (2k)(2k-1)Q_2(2k) z^{2k-2} \nonumber \\
& & =\sum_{k>2} (2k)^2 Q_2(2k) z^{2k-2} -\sum_{k>2} (2k)Q_2(2k)
z^{2k-2} \nonumber \\
& & =\sum_{k>2} (2k)^2 Q_2(2k) z^{2k-2} +O\left(\frac{1}{(1-z)^3
\log^2{\frac{1}{1-z}}}\right).
\end{eqnarray}
Combining (\ref{secder}), (\ref{mainterm}) and (\ref{asexpect}),
we obtain
$$
\sum_{k>2} (2k)^2 Q_2(2k) z^{2k-2} \sim\frac{6}{(1-z)^4
\log^2{\frac{1}{1-z}}}, \quad z\to 1^-.
$$
Applying again Hardy-Littlewood-Karamata Tauberian  theorem with
$r=1, \rho=4, L(t)=\frac{6}{\log^2{t}}$, we observe that
\begin{equation}\label{secmoment}
\sum_{2<k\le n} (2k)^2Q_2(2k)\sim\frac{4n^4}{\log^2{n}}, \quad
n\to\infty.
\end{equation}
Similarly to (\ref{expect}), for the second moment of $G_n$ we
have
\begin{equation}\label{gnsecmoment}
\frac{1}{n^2} \mathbb{E}(G_n^2) =\sum_{2<k\le n}
\left(\frac{k}{n}\right)^2 \mathbb{P}(G_n=k) =\sum_{2<k\le n}
\left(\frac{2k}{2n}\right)^2 \frac{Q_2(2k)}{\sum_{2<k\le n}
Q_2(2k)}.
\end{equation}
Now, from (\ref{secmoment}), (\ref{gnsecmoment}), (\ref{cardinal})
and (\ref{asympcard}) it follows that
$$
\lim_{n\to\infty}\frac{1}{n^2}\mathbb{E}(G_n^2)=\frac{1}{2},
$$
which, in combination with the result of Theorem 1(i), completes
the proof of part (ii). $\rule{2mm}{2mm}$

\end{document}